\documentclass[12pt]{article}
  \usepackage{amsfonts}
  \usepackage[square,sort&compress,comma,numbers]{natbib} 
   \def\R{\mathbb{R}}
   \def\N{\mathbb{N}}
   
   \def\1{{\rm I\mskip -10.5mu 1}} 
  
   \def\io{\int_\Omega}

   \def\H{{H^1_0(\Omega)}}

   \def\b{{\beta}}
   \def\e{{\varepsilon}}
   \def\be{{\beta,\e}}

   \def\supp{\mathop{\rm supp}\nolimits}

   \def\diam{\mathop{\rm diam}\nolimits}

   \def\lg{\mathop{\rm lg}\nolimits}
   \def\cp{\mathop{\rm cap}\nolimits}
 
   \def\no{\noindent}
   \def\proof{\mbox {{\underline {\sf Proof}} \hspace{2mm}}}
   \def\qed{{\hfill {\em q.e.d.}\\\vspace{1mm}}}
   
   \newcommand{\beq}{\begin{equation}}
   \newcommand{\eeq}{\end{equation}}

\newtheorem{df}{Definition}[section]
\newtheorem{prop}[df]{Proposition}

\newtheorem{teo}[df]{Theorem}
\newtheorem{rem}[df]{Remark}

\newtheorem{cor}[df]{Corollary}

 \newcommand{\sezione}[1]{\section{#1}\setcounter{equation}{0}}

\begin{document}


   \title{Variational properties of the first curve of the
     Fu\v{c}\'{\i}k spectrum for elliptic operators}

 \date{}

  \maketitle


\vspace{-2cm}

\begin{center}

{ {\bf Riccardo MOLLE$^a$,\quad Donato PASSASEO$^b$}}

\vspace{5mm}

{\em
${\phantom{1}}^a$Dipartimento di Matematica,
Universit\`a di Roma ``Tor Vergata'',\linebreak
Via della Ricerca Scientifica n. 1,
00133 Roma, Italy.\linebreak
{\tt  molle@mat.uniroma2.it  }}

\vspace{2mm}

{\em
${\phantom{1}}^b$Dipartimento di Matematica ``E. De Giorgi'',
  Universit\`a di Lecce,\linebreak 
P.O. Box 193, 73100 Lecce, Italy.
}
\end{center}

\vspace{5mm}


{\small {\sc \noindent \ \ Abstract.} - 
In this paper we present a new variational characteriztion of the
first nontrival curve of the Fu\v{c}\'{\i}k spectrum for elliptic
operators with Dirichlet boundary conditions.
Moreover, we describe the asymptotic behaviour and some properties of
this curve and of the corresponding eigenfunctions.
In particular, this new characterization allows us to compare the
first curve of the  Fu\v{c}\'{\i}k spectrum with the infinitely many
curves we obtained in previous works (see \cite{c,im}): for example,
we show that these curves are all asymptotic to the same lines as the
first curve, but they are all distinct from such a curve.
\footnote{ POST-PRINT: 
https://link.springer.com/article/10.1007/s00526-015-0920-4  }
\vspace{2mm}



\vspace{1mm}



}


\sezione{Introduction}


The  Fu\v{c}\'{\i}k spectrum, first introduced in \cite{F2} and
\cite{D77}, plays an important role in the study of some elliptic
problems with linear growth.
Let us consider, for example, the Dirichlet problem
\beq
\label{D}
\Delta u+g(x,u)=0\quad\mbox{ in }\Omega,\qquad u=0\quad\mbox{ on
}\partial\Omega,
\eeq
where $\Omega$ is a bounded connected domain of $\R^N$ with $N\ge 1$
and $g$ is a Carath\'eodory function in $\Omega\times\R$ such that
\beq
\lim_{t\to -\infty}{g(x,t)\over t}\, =\alpha,\qquad\lim_{t\to
  +\infty}{g(x,t)\over t}\, =\beta\qquad\forall x\in\Omega
\eeq
with $\alpha,\beta$ in $\R$.
Existence and multiplicity of solutions for problems of this type are
strictly related to the position of the pair
$(\alpha,\beta)$ with respect to the  Fu\v{c}\'{\i}k
spectrum $\Sigma$ which is defined as the set of all the pairs
$(\alpha,\beta)\in\R^2$ such that the Dirichlet problem 
\beq
\label{F}
\Delta u-\alpha u^-+\beta u^+=0\quad\mbox{ in }\Omega, \qquad
u=0\quad\mbox{ on }\partial \Omega
\eeq
has nontrivial solutions (here, $u^\pm=\max\{\pm u,0\}$ and $u$
nontrivial means $u\in H^1_0(\Omega)$, $u\not\equiv 0$).
In fact, these problems may lack compactness in the sense that the
well known Palais-Smale compactness condition fails if the pair
$(\alpha,\beta)$ belongs to the  Fu\v{c}\'{\i}k spectrum $\Sigma$;
moreover, the topological properties of the sublevels of the
corresponding energy functional depend on the position of
$(\alpha,\beta)$ with respect to $\Sigma$.

\no After the pionering papers \cite{C} and \cite{AP} on these
problems, the important role of the  Fu\v{c}\'{\i}k spectrum has been
pointed out in \cite{F2} and \cite{D77}.
Then, several works have been devoted to describe the structure of
$\Sigma$ (see, for example, 
\cite{AC,BFS,Be,Ca,Ca2,CG,D77,D81,D5,dFG,F,F2,FK,GK,GK82,GM,HR,LLLP,Ma,MM,R,Ru,Sc1,Sc2,S2} 
etc. \ldots).

\no Let us denote by $\lambda_1<\lambda_2\le\lambda_3\le\ldots$ the
eigenvalues of $-\Delta$ in $H^1_0(\Omega)$.
It is clear that $\Sigma$ includes the lines $\{\lambda_1\}\times\R$
and $\R\times\{\lambda_1\}$, contains all the pairs
$(\lambda_i,\lambda_i)$ $\forall i\in \N$ (that are the only pairs
$(\alpha,\beta)$ of $\Sigma$ such that $\alpha=\beta$) and is
symmetric with respect to the line $\{(\alpha,\beta)\in\R^2\ :\
\alpha=\beta\}$ (because a function $u$ satisfies (\ref{F}) if and
only if $-u$  satisfies (\ref{F}) with $(\beta,\alpha)$ in place of
$(\alpha,\beta)$).
Moreover, if $\alpha\neq\lambda_1$, $\beta\neq\lambda_1$ and
$(\alpha,\beta)\in\Sigma$, then $\alpha>\lambda_1$, $\beta>\lambda_1$
and the eigenfunctions corresponding to the pair $(\alpha,\beta)$ are
sign changing functions. 

\no In the case $N=1$, $\Sigma$ is completely known and may be obtained by
direct computation.
It consists of curves emanating from the pairs
$(\lambda_i,\lambda_i)$; if $i$ is an even positive integer, there
exists only one curve while, if $i$ is odd, there exist exactly two
curves emanating from $(\lambda_i,\lambda_i)$.
All these curves are smooth, unbounded and decreasing.
Moreover, on each curve, $\alpha$ tends to an eigenvalue of $-\Delta$
in $H^1_0(\Omega)$ as $\beta\to+\infty$.
Conversely, for every eigenvalue $\lambda_i$ there exist exactly three
curves asymptotic to the lines $\{\lambda_i\}\times\R$ and
$\R\times\{\lambda_i\}$; they pass, respectively, through the pairs
$(\lambda_{2i-1},\lambda_{2i-1})$, $(\lambda_{2i},\lambda_{2i})$
and $(\lambda_{2i+1},\lambda_{2i+1})$.

\no On the contrary, in the case of higher dimensions the known
results are much less complete and the description of $\Sigma$ remains
a largely open question.
It is known that $\Sigma$ is a closed set and that the lines
$\{\lambda_1\}\times \R$ and $\R\times \{\lambda_1\}$ (the trivial
part of $\Sigma$) are isolated in $\Sigma$ (see \cite{D77}).
Many results concern the curves of $\Sigma$ emanating from each pair
$(\lambda_i,\lambda_i)$ (local existence and multiplicity, local and
global properties, etc. \ldots).
In particular, if $\lambda_i$ has multiplicity $m$,
i.e. $\lambda_{i-1}<\lambda_i=\ldots=\lambda_{i+m-1}<\lambda_{i+m}$,
from the point $(\lambda_i,\lambda_i)$ arises a continuum composed by
a lower and an upper curve, both decreasing (and maybe coincident)
such that $\Sigma\cap(\lambda_{i-1},\lambda_{i+m})^2$ lies between
these two curves, so the open squares $(\lambda_{i-1},\lambda_{i})^2$ 
and $(\lambda_{i},\lambda_{i+m})^2$ do not contain any point of
$\Sigma$ (see, for example, \cite{GK,Ma,Sc1,Sc2} and the references
therein). 

\no Combining these results, one can infer that $\Sigma$ contains a
first nontrivial curve, which passes through $(\lambda_2,\lambda_2)$
and extends to infinity.
In \cite{dFG} the authors prove directly the existence  of such a
first curve, give a variational characterization of it and derive
several properties.
In particular, they show that this curve is asymptotic to the lines 
$\{\lambda_1\}\times \R$ and $\R\times \{\lambda_1\}$, give a new
proof of the fact that these lines are isolated in $\Sigma$ and deduce
that all the eigenfunctions corresponding to points of the first curve
have exactly two nodal domains (extending the well known Courant nodal
domains theorem).

\no Recently (see \cite{J1,J2,J3}) we have obtained new existence and
multiplicity results for a class of Dirichlet problems of type
(\ref{D}) (in particular for semilinear problems with jumping
nonlinearities) using a variational method that does not require to
know whether or not the pair $(\alpha,\beta)$ belongs to $\Sigma$ and,
in addition, may be used to give new information on the structure of
$\Sigma$.
In fact (see \cite{c,im}) using this method we have proved that if
$N\ge 2$ there exist infinitely many curves in $\Sigma$, asymptotic to the lines
$\{\lambda_1\}\times \R$ and $\R\times \{\lambda_1\}$ 
(while, if $N=1$, $\Sigma$ has only two curves asymptotic to these lines).
More precisely (see also Theorem \ref{T3.1}) we have proved that, if
$N\ge 2$ and $k\in \N$, for $\beta>0$ large enough there exists
$\alpha_{k,\beta}>\lambda_1$ such that
$(\alpha_{k,\beta},\beta)\in\Sigma$; moreover, for all $k\in\N$,
$\lim_{\beta\to +\infty} \alpha_{k,\beta}=\lambda_1$,
$\alpha_{k,\beta}$ depends continuously on $\beta$ and
$\alpha_{k,\beta}<\alpha_{k+1,\beta}$ (notice that the method
developed in \cite{J1,J2,J3} has been also used in the study of some
nonlinear scalar field equations (see \cite{CPS1,CPS2,CMP})).

\no The following natural question remains still open: where do these
curves come from? 
Most probably, they come from bifurcations of the first nontrivial
curve of $\Sigma$ or from pairs $(\lambda_i,\lambda_i)$ of higher
eigenvalues.

\no The results obtained in \cite{MR} seem to confirm our conjecture.
In fact, in \cite{MR} the authors study the  Fu\v{c}\'{\i}k spectrum
of the laplacian on a two-dimensional torus $T^2$ with periodic
conditions and, exploiting the invariance properties of $T^2$ with respect to
translations, they prove that at least two global curves emanate from
every pair of eigenvalues: a global curve which can be obtained
explicitely and a global curve which can be characterized
variationally using a suitable topological index (see
\cite{B,BLMR,Mar}).
The explicit curves are asymptotic to  lines $\{\lambda\}\times \R$
and $\R\times \{\lambda\}$ for suitable constants $\lambda>0$, while
the variational curves are all asymptotic to $\{0\}\times \R$
and $\R\times \{0\}$, the trivial lines of the  Fu\v{c}\'{\i}k
spectrum in $T^2$.
Therefore, the explicit and the variational curves cannot coincide
globally and many curve crossings must occur.
Moreover, on the first explicit curve there exist infinitely many
points of bifurcations (associated to symmetry breaking phenomena).

\no These results suggest that, in higher dimensions, $\Sigma$ has a
quite complicated structure even if, for example, $\Omega$ is a ball
of $\R^N$ with $N\ge 2$ (where we observed similar phenomena that
motivated our conjecture).

\no In the present paper we give a new variational characterization of
the first nontrivial curve of $\Sigma$.
We prove that the infimum
\beq
\inf\left\{\int_\Omega |Du^-|^2dx\ : \ u\in H^1_0(\Omega),\
\|u^+\|_{L^2(\Omega)}=\|u^-\|_{L^2(\Omega)}=1,
 \int_\Omega|Du^+|^2dx =\beta\right\}
\label{inf}
\eeq
is achieved for all $\beta>\lambda_1$; moreover we prove that, if we denote by
$\alpha_\beta$ the above infimum, then $\alpha_\beta$ is continuous
and strictly decreasing with respect to $\beta$ in
$]\lambda_1,+\infty[$, $\lim_{\beta\to +\infty}\alpha_\beta=\lambda_1$,
$\lim_{\beta\to \lambda_1}\alpha_\beta=+\infty$,
$(\alpha_\beta,\beta)\in\Sigma$ $\forall\beta>\lambda_1$ and 
\beq
\{(\alpha,\beta)\in\Sigma\ :\ \alpha>\lambda_1,\ \beta>\lambda_1\}
\subseteq 
\{(\alpha,\beta)\in\R^2\ :\ \beta>\lambda_1,\ \alpha\ge \alpha_\beta\
\forall\beta>\lambda_1\}.
\eeq
Therefore (using $\beta$ as parameter) the first nontrivial curve of
$\Sigma$ may be described as the set $\{(\alpha,\beta)\in\R^2$ : $\b>\lambda_1$,
$\alpha=\alpha_\beta$ $\forall\beta>\lambda_1\}$ (a similar
description holds if we use $\alpha$ as parameter). 

\no The eigenfunctions corresponding to $(\alpha_\beta,\beta)$ are
suitable smooth functions $u_\beta\in H^1_0(\Omega)$, such that
$u_\beta^+\not\equiv 0$, $u_\beta^-\not\equiv 0$ in $\Omega$ and the
function $\bar
u_\beta=-\|u_\beta^-\|^{-1}_{L^2(\Omega)}u_\beta^-+\|u_\beta^+\|^{-1}_{L^2(\Omega)}u_\beta^+$
is a minimizing function for (\ref{inf}).

\no This variational characterization of the first nontrivial curve of
$\Sigma$, which is different from the ones obtained in \cite{dFG} and
in \cite{MR}, has been first announced in \cite{l}.

\no All the properties of the first curve and of the corresponding
eigenfunctions may be easily deduced from this new characterization. 
In the present paper, in particular, we deduce that all the curves
we obtained in \cite{im} are distinct from the first curve.
In fact, suitable asymptotic estimates of $\alpha_\beta$, as $\beta\to
+\infty$, allow us to say that for every positive integer $k$ we have
$\alpha_\beta<\alpha_{k,\beta}$ when $\beta$ is large enough.

\no The asymptotic behaviour of the eigenfunctions $u_\beta$
corresponding to $(\alpha_\beta,\beta)$, as $\beta\to +\infty$, shows
that the support of $u_\beta^+$ is localized near the boundary of
$\Omega$ for $\beta$ large enough and that
$\|u_\beta\|^{-1}_{L^2(\Omega)} u_\beta\to -e_1$ in $\Omega$, where
$e_1$ denotes the positive eigenfunction of $-\Delta$ in
$H^1_0(\Omega)$, corresponding to $\lambda_1$ and normalized in
$L^2(\Omega)$.
On the contrary, the eigenfunctions $u_{k,\beta}$ corresponding to the
pairs $(\alpha_{k,\beta},\beta)$, for $\beta$ large enough, have the
support of $u_{k,\beta}^+$  localized near the maximum points of $e_1$
(see \cite{c,im} and also Theorem \ref{T3.1} and Proposition \ref{P3.2}).
This fact suggests that, arguing as in \cite{im}, it is possible to
construct a new class of infinitely many curves of $\Sigma$,
asymptotic to the lines $\{\lambda_1\}\times \R$ and $\R\times
\{\lambda_1\}$, corresponding to eigenfunctions having an arbitrarily
large number of bumps localized near the boundary of $\Omega$ (see
Remark \ref{R3.5} for more details about this construction).

\no The method we use in this paper is completely variational. 
For all $\beta>\lambda_1$, the eigenfunction $u_\beta$ corresponding
to the pair $(\alpha_\beta,\beta)$ is obtained as critical point of
the functional $f(u)=\int_\Omega[|Du|^2-\beta(u^+)^2]dx$ constrained
on the set $S=\{u\in H^1_0(\Omega)$ : $\int_\Omega (u^-)^2dx=1\}$
(here, $\alpha_\beta$ arises as the Lagrange multiplier with respect
to the constraint $S$).
In Section 2 we introduce also another functional $f_{\beta,\e}$,
converging to $f_\beta$ as $\e\to 0$, which for all $\e>0$ presents
more manageable variational properties with respect to $f_\beta$;
thus we first find constrained critical points for $f_{\b,\e}$ and
then we let $\e\to 0$ in order to obtain the variational characterization of
the first nontrivial curve of $\Sigma$ (see Theorem \ref{T2.1}).
In Section 3 we analyse the asymptotic behaviour, as $\beta\to
+\infty$, of this curve and of the corresponding eigenfunctions, we
compare this curve with the infinitely many curves obtained in
\cite{im} and we discuss some possible generalizations, forthcoming
results on related questions, etc. \ldots


\sezione{Variational characterization of the first curve of $\Sigma$}


The following theorem gives a variational characterization of the
first nontrivial curve of the  Fu\v{c}\'{\i}k spectrum $\Sigma$.

\begin{teo}
\label{T2.1}
Let $\Omega$ be a bounded connected domain of $\R^N$ with $N\ge 1$.
For all $\b>\lambda_1$, let us set
\beq
\label{ab}
\alpha_\b=\inf\left\{ \int_\Omega |Du^-|^2dx\, :\, u\in H^1_0(\Omega),\
  \|u^+\|_{L^2(\Omega)}=\|u^-\|_{L^2(\Omega)}=1, \int_\Omega |Du^+|^2dx=\b\right\}.
\eeq
Then, $\alpha_\beta>\lambda_1$, $(\alpha_\b,\b)\in\Sigma$
$\forall\b>\lambda_1$ and $\alpha_\b\le\alpha$ for every
$\alpha>\lambda_1$ such that $(\alpha,\b)\in\Sigma$.
Moreover, $\alpha_\b$ is continuous and strictly decreasing with
respect to $\beta$ in $]\lambda_1,+\infty[$,
$\alpha_{\lambda_2}=\lambda_2$, the infimum in (\ref{ab}) is achieved
$\forall \b>\lambda_1$ and an eigenfunction corresponding to the pair
$(\alpha_\b,\b)$ is given by $u_\b=-\bar u_\b^-+\mu_\b\bar u_\b^+$,
where $\bar u_\b$ is a minimizing function for (\ref{ab}) and $\mu_\b$
is a suitable positive constant.

\no As $\b\to +\infty$, $\alpha_\b\to\lambda_1$ and $u_\b^-\to e_1$ in
$H^1_0(\Omega)$; as $\b\to\lambda_1$, $\alpha_\b\to +\infty$
and $\|u^+_\beta\|^{-1}_{L^2(\Omega)}u^+_\b\to e_1$ in $H^1_0(\Omega)$.
\end{teo}

In order to prove this theorem, for all $\b>0$ and $\e>0$ we consider
the functional $f_{\b,\e}:H^1_0(\Omega)\to\R$ defined by
\beq
\label{fbe}
f_{\b,\e}(u)=\int_\Omega |D u|^2dx-2\int_\Omega G_{\b,\e}(u)dx,
\eeq
where $G_{\b,\e}(t)=\int_0^t g_{\b,\e}(\tau)d\tau$ $\forall t\in\R$,
with $g_{\b,\e}(\tau)=0$ $\forall \tau\le \e$ and
  $g_{\b,\e}(\tau)=\beta(\tau-\e)$ $\forall\tau\ge \e$.
Our first aim is to find sign changing functions $u\in H^1_0(\Omega)$
which are constrained critical points for the functional $f_{\b,\e}$
constrained on the set $S=\{u\in H^1_0(\Omega)$ : $ \int_\Omega
(u^-)^2dx=1\}$.
Therefore, we consider the set $M^{\b,\e}$ consisting of all the
functions $u$ in $S$ such that $u^+\not\equiv 0$ and
$f'_{\b,\e}(u)[u^+]=0$.

\no One can easily verify that for all $\e>0$, if a sign  changing 
function $u\in S$ is a critical point for $f_{\b,\e}$ constrained on
$S$, then $u\in M^{\b,\e}$ and $f'_{\b,\e}(-u^-+tu^+)[u^+]$ is
positive for $t\in]0,1[$ and negative for $t>1$ (because
${1\over\tau}\, g(\tau)$ is strictly increasing with respect to $\tau$
in $]\e,+\infty[$); so the function $u$ is the unique maximum point
for $f_{\b,\e}$ on the set $\{-u^-+tu^+$ : $t\ge 0\}$ which is
included in $S$.

\no Notice that, for $\e=0$, the functional $f_\b=f_{\b,0}$ and the
sets $M^{\b,0}$ do not have the same properties; this is the reason for
introducing first the functional $f_{\b,\e}$ which, for all
$\e>0$, presents more manageable variational properties and then we
let $\e\to 0$.

\begin{prop}
\label{P2.2}
For all $\e>0$, we have $M^{\b,\e}\neq\emptyset$ if and only if
$\b>\lambda_1$.
If $\b>\lambda_1$, the minimum of $f_{\b,\e}$ on $M^{\b,\e}$ is
achieved for all $\e>0$ and every minimizing function $u_{\b,\e}$
satisfies the equation
\beq
\label{eq3}
\Delta u_{\b,\e}+g_{\b,\e}(u_{\b,\e})-\alpha_{\b,\e}u^-_\be=0\qquad\mbox{
  in }\Omega,
\eeq
where $\alpha_{\b,\e}=\int_\Omega|Du_{\b,\e}^-|^2 dx$ $\forall\e>0$,
$\forall \b>\lambda_1$.
\end{prop}

\proof
Notice that $g_{\b,\e}(\tau)\le\b\tau$ $\forall\tau\ge 0$,
$\forall\b > 0$, $\forall \e >0$; moreover, if a function $u\in S$
satisfies $\io (u^+)^2dx > 0$, then  $\io |Du^+|^2dx>\lambda_1\io
(u^+)^2dx$.
Thus, if $\b\le\lambda_1$ we obtain $\io |Du^+|^2dx>\b\io(u^+)^2dx\ge\io
g_{\b,\e}(u^+)u^+dx$, namely $f'_{\b,\e}(u)[u^+]>0$, which implies
$u\not\in M^{\b,\e}$.
Therefore, if $\b\le\lambda_1$, we have $M^{\b,\e}=\emptyset$
$\forall\e>0$.

\no On the contrary, if $\b>\lambda_1$ we have
$M^{\b,\e}\neq\emptyset$. 
In fact, if $\b>\lambda_1$, one can easily construct a function $\bar
u\in S$ such that $\io |D\bar u^+|^2 dx<\b\io (\bar u^+)^2dx$.
For example, fix $x_0\in \partial\Omega$, set
$\Omega_r=\Omega\setminus\overline B(x_0,r)$,
$\tilde\Omega_r=\Omega\cap B(x_0,r)$ and, for $r\in
]0,(1/2)\diam(\Omega)[$, consider the positive eigenfunctions
$e_1(\Omega_r)$ and $e_1(\tilde\Omega_r)$, normalized in $L^2$,
corresponding to the first eigenvalues of the operator $-\Delta$ in
$H^1_0(\Omega_r)$ and $H^1_0(\tilde\Omega_r)$ respectively; then, the
function $\bar u$ such that $\bar u=e_1(\Omega_r)$ in $\Omega_r$ and 
$\bar u=-e_1(\tilde \Omega_r)$ in $\tilde\Omega_r$ has the required
properties for $r$ small enough.

\no Notice that $f_{\b,\e}(-\bar u^-+t\bar u^+)=f_{\b,\e}(-\bar
u^-)+f_{\b,\e}(t\bar u^+)$ $\forall t\ge 0$; moreover, one can easily
verify by direct computation that
\beq
\lim_{t\to +\infty}{1\over t^2} f_{\b,\e}(t\bar u^+)=\io |D\bar
u^+|^2dx-\b\io(\bar u^+)^2dx<0
\eeq
and (since $\e>0$)
\beq
\lim_{t\to 0}{1\over t^2} f_{\b,\e}(t\bar u^+)=\io |D\bar
u^+|^2dx>0.
\eeq
Therefore, we infer that for all $\e>0$ there exists $\bar t>0$ such
that $f_{\b,\e}(-\bar u^-+\bar t\bar u^+)\ge f_{\b,\e}(-\bar u^-+
t\bar u^+)$ $\forall t\ge 0$, which implies $-\bar u^-+\bar t\bar
u^+\in M^{\b,\e}$.
Thus, $M^{\b,\e}\neq\emptyset$ $\forall\b>\lambda_1$, $\forall\e>0$.

\no Now, let us prove that the infimum $\inf_{M^{\b,\e}}f_{\b,\e}$ is
achieved for all $\b>\lambda_1$ and $\e>0$.
Let us consider a minimizing sequence $(u_n)_n$.
Notice that $f_{\b,\e}(u_n)=f_{\b,\e}(- u_n^-)+f_{\b,\e}(u_n^+)$
$\forall n\in\N$, where $=f_{\b,\e}(-
u_n^-)=\io|Du_n^-|^2dx\ge\lambda_1$ (because
$\|u_n^-\|_{L^2(\Omega)}=1$) and $f_{\b,\e}(u_n^+)>0$ since $u_n\in
M^{\b,\e}$ implies $f_{\b,\e}(u_n^+)=\max\{f_{\b,\e}(tu_n^+)$ : $t\ge
0\}>0$ $\forall \e>0$ (because ${1\over\tau}g_{\b,\e}(\tau)$ is
strictly increasing with respect to $\tau$ in $]\e,+\infty[$).
Taking into account that $\sup\{f_{\b,\e}(u_n)$ : $n\in\N\}<+\infty$,
it follows that 
\beq
\label{a}
\lambda_1\le\liminf_{n\to \infty}
f_{\b,\e}(-u_n^-)\le\limsup_{n\to\infty}f_{\b,\e}(-u_n^-)<+\infty
\eeq
and
\beq
\label{b}
0\le\liminf_{n\to \infty}
f_{\b,\e}(u_n^+)\le\limsup_{n\to\infty}f_{\b,\e}(u_n^+)<+\infty.
\eeq
Since $f_{\b,\e}(-u_n^-)=\io |Du^-_n|^2dx$, (\ref{a}) implies that the
sequence $(u_n^-)_n$ is bounded in $\H$.
Now, let us prove that also the sequence $(u_n^+)_n$ is bounded in
$\H$.
Taking into account that $f'_{\b,\e}(u_n)[u_n^+]=0$ $\forall n\in\N$,
we have
\beq
\label{f'}
\io |Du_n^+|^2dx =\io g_{\b,\e}(u^+_n)u_n^+dx\le\b\io (u^+_n)^2dx.
\eeq
Therefore, it suffices to prove that the sequence $(u_n^+)_n$ is
bounded in $L^2(\Omega)$.
Arguing by contradiction, assume that (up to a subsequence)
$\lim_{n\to\infty}\io (u_n^+)^2dx=\infty$ and set
$v_n=\|u_n^+\|^{-1}_{L^2(\Omega)}u_n^+$.
Then, (\ref{f'}) implies $\io |Dv_n|^2dx\le\b$ $\forall n\in\N$.
So (up to a subsequence) $(v_n)_n$ converges in $L^2(\Omega)$, weakly
in $\H$ and a.e. in $\Omega$ to a function $v\in\H$.
It follows that $\io |Dv|^2dx\le\b$, $\io v^2dx=1$ and $v\ge 0$ in
$\Omega$.
Moreover, taking into account that $f'_{\b,\e}(u_n)[u_n^+]=0$ $\forall
n\in\N$, and that $\lim_{n\to\infty}\io(u_n^+)^2dx=\infty$, one can
verify by direct computation that $\lim_{n\to\infty}\io |Dv_n|^2dx
=\b$.
It follows that
\beq
\lim_{n\to\infty}f'_{\b,\e}(tv_n)[v_n]=2t\b-2\io
g_{\b,\e}(tv)vdx\qquad\forall t\ge 0.
\eeq
Since $\io v^2dx=1$, for all $\e>0$ we obtain
\beq
\label{>}
\liminf_{t\to+\infty}\left[t\beta-\io
g_{\b,\e}(tv)v\,dx\right]=\liminf_{t\to+\infty}\io[\b tv-g_{\b,\e}(tv)]v\,dx>0.
\eeq
Now, let us set $t_n=\|u_n^+\|_{L^2(\Omega)}$ and notice that
$f'_{\b,\e}(tv_n)[v_n]>0$ $\forall t\in ]0,t_n[$ (because
${1\over\tau}g_{\b,\e}(\tau)$ is strictly increasing with respect to
$\tau$ in $]\e,+\infty[$).
Since we are assuming $\lim_{n\to\infty}t_n=+\infty$, we obtain
\beq
\liminf_{n\to\infty}f_{\b,\e}(u^+_n)=\liminf_{n\to\infty}\int_0^{t_n}f'_{\b,\e}(tv_n)[v_n]dt\ge
2\int_0^\tau\left[t\b-\io g_{\b,\e}(tv)v\,dx\right]dt\quad\forall\tau>0.
\eeq
Then, as $\tau\to+\infty$, from (\ref{>}) we obtain
$\lim_{n\to\infty}f_{\b,\e}(u^+_n)=+\infty$, in contradiction with
(\ref{b}).
Therefore, we can say that also the sequence $(u_n^+)_n$ is bounded in
$\H$.
It follows that there exists $u\in \H$ such that (up to a subsequence)
$(u_n)_n$ converges to $u$ in $L^2(\Omega)$, weakly in $\H$ and
a.e. in $\Omega$.
As a consequence of the $L^2(\Omega)$ convergence, we have
$\io(u^-)^2dx=1$.
Let us prove that $u^+\not\equiv 0$.
Arguing by contradiction, assume that $u^+\equiv 0$.
Then (because of the $L^2(\Omega)$ convergence) from (\ref{f'}) we
infer that $\lim_{n\to\infty}\io|Du_n^+|^2dx=0$, so we have
$\lim_{n\to\infty}f_{\b,\e}(u_n^+)=0$.
Therefore, we obtain a contradiction if we prove that
\beq
\label{inf_1}
\inf\{f_{\b,\e}(w^+)\ :\ w\in M^{\b,\e}\}>0\qquad\forall\e>0.
\eeq
Since $w\in M^{\b,\e}$ implies $f_{\b,\e}(w^+)=\max\{f_{\b,\e}(tw^+)$
: $t>0\}$, it is clear that it suffices to prove that there exist two
positive constants $\rho_{\b,\e}$ and $c_{\b,\e}$ such that
$f_{\b,\e}(w)\ge c_{\b,\e}$  $\forall w\in\H$ such that $\io
|Dw|^2dx=\rho_{\b,\e}$.

\no In order to prove the existence of these constants $c_{\b,\e}$ and
$\rho_{\b,\e}$, notice that, since $\b>\lambda_1$, there exists
$j_\b\in\N$ such that $\lambda_{j_\b}\le\beta<\lambda_{j_\b+1}$.
Let us denote by $S^1_\b$ and $S^2_\b$ the closed subspaces of
$\H$ spanned by the eigenfunctions of the Laplace operator $-\Delta$
in $\H$, corresponding to eigenvalues $\lambda_j$ with, respectively,
$1\le j\le j_\b$ and $j\ge j_\b+1$.
For all $\b>\lambda_1$ and $\e>0$, there exists $\nu_{\b,\e}>0$ such
that, if $w\in S^1_\beta$ and $\io |Dw|^2dx\le\nu_{\b,\e}^2$, then
$|w(x)|\le \e$ $\forall x\in\Omega$.

\no For all $w\in\H$ such that $\io |Dw|^2dx\le\nu_{\b,\e}^2$, set
$w=w_{1,\b}+w_{2,\b}$, with $w_{1,\beta}\in S^1_\beta$ and
$w_{2,\beta}\in S^2_\b$.
Then, taking into account that $\io |Dw_{1,\b}|^2dx \le \nu^2_{\b,\e}$
and as a consequence $w_{1,\b}\le\e$ in $\Omega$, we have
$f_{\b,\e}(w_{1,\b})=\io |Dw_{1,\b}|^2dx$ and
$f'_{\b,\e}(w_{1,\beta})[w_{2,\b}]=0$.
Therefore, we obtain
\beq
f_{\b,\e}(w)=f_{\b,\e}(w_{1,\b}+w_{2,\b})=f_{\b,\e}(w_{1,\b}+w_{2,\b})-f_{\b,\e}(w_{1,\b})+\io
|Dw_{1,\b}|^2dx,
\eeq
where
\begin{eqnarray}
f_{\b,\e}(w_{1,\b}+w_{2,\b})-f_{\b,\e}(w_{1,\b}) 
& \ge & 
f'_{\b,\e}(w_{1,\b})[w_{2,\b}]+\io |Dw_{2,\b}|^2dx-\b \io w_{2,\b}^2dx
\nonumber\\
&=&
\io [|Dw_{2,\b}|^2-\b w_{2,\b}^2]dx
\nonumber\\
&\ge&
\left( 1-{\b\over \lambda_{j_\b+1}}\right) \io |Dw_{2,\b}|^2dx
\end{eqnarray}
because $\io |Dw_{2,\b}|^2dx\ge \lambda_{j_\b+1} \io w_{2,\b}^2dx$.
It follows that, for a suitable $\tilde c_{\b,\e}>0$, we have
$f_{\b,\e}(w)\ge\tilde c_{\b,\e}\io |Dw|^2dx$ $\forall w\in\H$ such
that $\io |Dw|^2dx\le\nu^2_{\b,\e}$.
Therefore, it follows easily that there exist two constants
$\rho_{\b,\e}\in]0,\nu_{\b,\e}[$ and $c_{\b,\e}>0$ satisfying the
required properties.
Thus, we can say that $u^+\not\equiv 0$.

\no From the weak $\H$ convergence, it follows that
$f_{\b,\e}'(u^+)[u^+]\le 0$; on the other hand, a direct computation
shows that $\lim_{t\to 0}{1\over t} f'_{\b,\e}(tu^+)[u^+]=2\io
|Du^+|^2dx$, so we infer that $f'_{\b,\e}(tu^+)[u^+]>0$ for $t>0$
small enough (because $u^+\not\equiv 0$).
Therefore, there exists $\tilde t\in]0,1]$ such that the function
$\tilde u=-u^-+\tilde t u^+$ belongs to $M^{\b,\e}$.
Moreover, since $f_{\b,\e}(\tilde t u^+_n)\le f_{\b,\e}(u^+_n)$ $\forall
n\in \N$, we have
\beq
\liminf_{n\to\infty} f_{\b,\e}(\tilde t u_n^+)\le 
\liminf_{n\to\infty} f_{\b,\e}( u_n^+).
\eeq
It follows that
\beq
f_{\b,\e}(\tilde u)\le \liminf_{n\to\infty}f_\be (-u_n^-+\tilde tu_n^+)\le
\liminf_{n\to\infty} f_{\b,\e}( u_n)
=\inf\{f_\be(u)\ :\ u\in M^\be\}.
\eeq
Thus, we can say that the infimum of $f_\be$ on $M^\be$ is achieved
and that $f_\be(\tilde u)=\min_{M^\be}f_\be$.

\no Let $\b>\lambda_1$, $\e>0$ and $u_\be$ be a minimizing function for
$f_\be$ on $M^\be$.
Our aim is to prove that $u_\be$ is a constrained critical point for
the functional $f_\be$ constrained on the set $S$, namely that there
exists a Lagrange multiplier $\alpha_\be$ such that 
\beq
\label{weak}
{1\over 2} f'_\be(u_\be)[\psi]=-\alpha_\be\io u^-_\be\psi\,
dx\qquad\forall\psi\in H^1_0(\Omega) 
\eeq
(that is, $u_\be$ is a weak solution of the equation (\ref{eq3})).

\no Let us point out that, unlike the case of the smooth constraint
$\io (u^-)^2dx=1$, for which the Lagrange multipliers theorem applies
(and gives the multiplier $\alpha_\be$), the constraint
$f'_\be(u)[u^+]=0$ does not satisfy the regularity conditions
  required in that theorem.
However, it is a ``natural constraint'' in the sense that it does not
give rise to Lagrange multipliers.

\no Notice that (as we observed before)
$f'_\be(u+tu^+)[u^+]$ is positive for $t\in ]-1,0[$ and negative for
$t>0$.
Therefore, $u_\be$ is the unique maximum point for $f_\be$ on the set
$\{u_\be+tu_\be^+$ : $t\ge -1\}$.
Then, arguing by contradiction, assume that (\ref{weak}) is not
satisfied for any choice of the multiplier $\alpha_\be$ in $\R$.
It follows by standard arguments that there exists a continuous map
$\eta:]-1,+\infty[\to\H$ such that $\eta(t)=u_\be+tu^+_\be$ if
$|t|>{1\over 2}$, $\|\eta(t)^-\|_{L^2(\Omega)}=1$,
$\eta(t)^+\not\equiv 0$, $f_\be(\eta(t))<f_\be(u_\be)$ $\forall t\ge
-1$.

\no Therefore, we infer that there exists $\bar t\in \left[-{1\over
    2},{1\over 2}\right]$ such that $\eta(\bar t)\in M^\be$, which
gives a contradiction because $f_\be (\eta(\bar t))<f_\be(u_\be)$ and 
$f_\be(u_\be)=\min_{M^\be}f_\be$.
Thus, we can conclude that there exists a multiplier $\alpha_\be$ in
$\R$ such that (\ref{weak}) holds.

\no Finally, notice that, if in (\ref{weak}) we set $\psi=u_\be^-$, we
easily obtain $\alpha_\be=\io |Du^-_\be|^2dx$, so the proof is
complete.

\qed

\no Now, our aim is to describe the behaviour as $\e\to 0$ of the
minimizing function $u_\be$ given by Proposition \ref{P2.2}.

\begin{prop}
\label{P2.3}
For all $\b>\lambda_1$ and $\e>0$, let $u_\be$ be a minimizing
function for the functional $f_\be$ on the set $M^\be$ and put $\bar
u_\be=-u^-_\be+\|u^+_\be\|_{L^2(\Omega)}^{-1}u^+_\be$.
Then, up to a subsequence, $\bar u_\be$ converges in $\H$, as $\e\to
0$, to a function $\bar u_\b$ such that $\io |D\bar u^+_\b|^2dx=\beta$
and
\begin{eqnarray}
\io |D\bar u^-_\b|^2dx&=&\min\left\{\io|Du^-|^2dx \ :\ u\in\H,\right.
\label{min} \\
& &\hspace{2cm} \left.
  \|u^+\|_{L^2(\Omega)}=\|u^-\|_{L^2(\Omega)}=1,\ \io|D u^+|^2dx=\b\right\}.
\nonumber
\end{eqnarray}
\end{prop}

\proof
Notice that, since $u_\be\in M^\be$, we have 
\beq
\io |Du^+_\be|^2dx=\io g_\be(u_\be)u^+_\be dx<\b\io
(u^+_\be)^2 dx\quad\forall\e>0.
\eeq
Therefore, we obtain $\io |D\bar u_\be^+|^2dx<\b$ $\forall\e>0$.

\no Moreover, we have also
\beq
\label{u-}
\limsup_{\e\to 0}\io |Du^-_\be|^2dx<+\infty.
\eeq
In fact, since $\b>\lambda_1$, there exists $\hat u\in \H$ such that
$\io(\hat u^+)^2dx=\io (\hat u^-)^2dx=1$, $\io |D\hat u^+|^2dx<\b$. 
As a consequence, we obtain
\beq
f_\be(u_\be)\le\io |D\hat u^-|^2dx+\max\{f_\be(t\hat u^+)\ :\ t\ge
0\}\qquad\forall\e>0
\eeq
and, as $\e\to 0$, since $\io |D\hat u^+|^2dx<\b$,
\beq
\label{hatu}
\limsup_{\e\to 0}f_\be(u_\be)\le\io |D\hat u^-|^2dx<+\infty,
\eeq
which implies (\ref{u-}) because $f_\be(u_\be)=\io
|Du^-_\be|^2dx+f_\be(u^+_\be)$ with $f_\be(u^+_\be)>0$ $\forall\e>0$.

\no It follows that, up to a subsequence, $\bar u_\be$ converges as
$\e\to 0$ to a function $\bar u_\b\in \H$ in $L^2(\Omega)$, weakly in
$\H$ and a.e. in $\Omega$.

\no Let us prove that, indeed, $\bar u_\be\to\bar u_\b$ strongly in
$\H$ as $\e\to 0$.
In fact, we have
\beq
\io |D\bar u_\b^+|^2dx=\b \quad\mbox{ and }\quad \lim_{\e\to 0}\io
|Du^-_\be|^2dx=\io |Du^-_\b|^2dx.
\eeq
For the proof, we argue by contradiction and assume that $\io |D\bar
u^+_\b|^2dx<\b$ or (up to a subsequence) $\io
|Du^-_\b|^2dx<\lim_{\e\to 0}\io |Du^-_\be|^2dx$.
In this case, by slight modifications of the supports of $u_\b^-$ and
$\bar u^+_\b$, one can construct a function $\tilde u_\b\in \H$ such
that $\|\tilde u_\b^-\|_{L^2(\Omega)}= \|\tilde
u_\b^+\|_{L^2(\Omega)}=1$, $\io |D\tilde u_\b^+|^2dx<\b$ and $\io
|D\tilde u_\b^-|^2dx<\lim_{\e\to 0}\io |Du_\be^-|^2dx$.

\no Then, for all $\e>0$, let us consider the function $\tilde
u_\be\in M^\be$ such that $\tilde u_\be^-=\tilde u_\b^-$ $\forall\e>0$
and $\tilde u^+_\be=t_\e\tilde u^+_\b$ where, for all $\e>0$, $t_\e$ is the
(unique) positive number such that $f'_\be(t_\e\tilde u^+_\b)[\tilde
u^+_\b]=0$ (such a number $t_\e$ exists because $\io |D\tilde
u^+_\b|^2dx<\b$). 

\no Thus, we have
\beq
f_\be(u_\be)-f_\be(\tilde u_\be)=f_\be(u^+_\be)-f_\be(\tilde
u^+_\be)+f_\be(u^-_\be)-f_\be(\tilde u^-_\be),
\eeq
where $f_\be(u^+_\be)\ge 0$ $\forall \e>0$, $\lim_{\e\to
  0}f_\be(\tilde u^+_\be)=0$ (because  $\io  |D\tilde u_\b^+|^2dx<\b$)
and $\lim_{\e\to 0}f_\be(u^-_\be)>\io |D\tilde
u_\b^-|^2dx=f_\be(\tilde u^-_\be)$ $\forall\e>0$.

\no It follows that $f_\be(u_\be)>f_\be(\tilde u_\be)$ for $\e>0$
small enough, which gives a contradiction because $\tilde u_\be\in
M^\be$ and $f_\be(u_\be)=\min_{M^\be}f_\be$.
Thus, we can conclude that $\bar u_\be\to\bar u_\b$ strongly in $\H$
as $\e\to 0$ and that $\io |D\bar u^+_\b|^2dx =\b$.

\no In a similar way, now we prove (\ref{min}). Arguing again by
contradiction, assume that there exists a function $v\in\H$ such that
$\|v^+\|_{L^2(\Omega)}=\|v^-\|_{L^2(\Omega)}=1$, $\io |Dv^+|^2dx=\b$
and $\io |Dv^-|^2dx<\io |D\bar u^-_\b|^2dx$.

\no In this case, by slight modifications of the supports of $v^+$ and
$v^-$, one can find $\hat v_\b\in\H$ such that $\|\hat
v^+_\b\|_{L^2(\Omega)}=\|\hat v_\b^-\|_{L^2(\Omega)}=1$, $\io |D\hat
v_\b^+|^2dx<\b$ and $\io |D\hat v^-_\b|^2dx<\io |D\bar u^-_\b|^2dx$.

\no Since $\io |D\hat v_\b^+|^2dx<\b$ and $\|\hat
v_\b^+\|_{L^2(\Omega)}=1$, it follows that for all $\e>0$ there exists
$\hat t_\e>0$ such that $f'_\be(\hat t_\e\hat v^+_\b)[\hat v^+_\b]=0$,
namely, the function $\hat v_\be=-\hat v_\b^-+\hat t_\e\hat v^+_\b$
belongs to $M^\be$.

\no Then, by direct computation, we obtain
\beq
f_\be(u_\be)-f_\be(\hat v_\be)=f_\be(u^-_\be)-f_\be(\hat
v^-_\be)+f_\be(u^+_\be)-f_\be(\hat v^+_\be),
\eeq
where $f_\be(u^+_\be)\ge 0$, $f_\be(\hat v^-_\be)=\io |D\hat
v^-_\b|^2dx$ $\forall\e>0$, $\lim_{\e\to 0} f_\be(\hat v^+_\be)=0$ and 
\beq
\lim_{\e\to 0}f_\be (u^-_\be)=\io |D\bar u_\b^-|^2dx>\io |D\hat
v^-_\b|^2dx.
\eeq
It follows that $f_\be(u_\be)>f_\be(\hat v_\be)$ for $\e>0$ small
enough; so we have again a contradiction because $\hat v_\be\in M^\be$
and $f_\be(u_\be)=\min_{M^\be}f_\be$ $\forall \e>0$.

\qed

\begin{prop}
\label{P2.4}
For all $\b>\lambda_1$ and $\e>0$, let $u_\be$ be a minimizing
function for the functional $f_\be$ on the set $M^\be$.
Then, as $\e\to 0$ (up to a subsequence) $u_\be$ converges in $\H$ to
a sign changing function $u_\b$ which solves the equation 
\beq
\Delta u_\b-\alpha_\b u^-_\b+\b u^+_\b=0\qquad\mbox{ in }\Omega,
\eeq
where $\alpha_\b$ is the positive number introduced in Theorem
\ref{T2.1}.
Moreover, the function $\bar
u_\b=-u_\b^-+\|u^+_\b\|_{L^2(\Omega)}^{-1}u_\b^+$ is a minimizing
function for $\alpha_\b$ (see (\ref{ab})) and
$\lambda_1<\alpha_\beta\le\alpha$ for every $\alpha>\lambda_1$ such
that $(\alpha,\beta)\in\Sigma$.
\end{prop}

\proof
As we proved in Proposition \ref{P2.2}, for all $\b>\lambda_1$ and
$\e>0$, $u_\be$ is a weak solution of the equation
\beq
\label{eq}
\Delta u_\be-\alpha_\be u^-_\be+g_\be(u_\be)=0\qquad\mbox{ in }\Omega,
\eeq
where $\alpha_\be=\io |Du^-_\be|^2dx$.

\no Moreover, by Proposition \ref{P2.3}, $\bar
u_\be=-u_\be^-+\|u^+_\be\|^{-1}_{L^2(\Omega)}u^+_\be$ converges, as
$\e\to 0$, to a function $\bar u_\b$ in $\H$.

\no Let us prove that
\beq
\liminf_{\e\to 0}\|u^+_\be\|_{L^2(\Omega)}>0\qquad\forall\b>\lambda_1.
\eeq
Arguing by contradiction, assume that (up to a subsequence)
$\lim_{\e\to 0}\|u^+_\be\|_{L^2(\Omega)}=0$ for some
$\beta>\lambda_1$.
In this case, since $f'_\be(u_\be)[u^+_\be]=0$, we infer that
$u^+_\be\to 0$ in $\H$ as $\e\to 0$.
Therefore, if we let $\e\to 0$, from (\ref{eq}) we obtain 
\beq
\io[D\bar u^-_\b D\psi-\bar\alpha_\b \bar
u_\beta^-\psi]dx=0\qquad\forall\psi\in\H,
\eeq
where $\bar\alpha_\b=\io |D\bar u_\b^-|^2dx$.
Thus, we have a contradiction because $D\bar u_\b^-\not\equiv 0$ on
$\Omega\cap\partial(\supp\bar u^-_\b)$.

\no Now, let us prove that
\beq
\limsup_{\e\to
  0}\|u^+_\be\|_{L^2(\Omega)}<+\infty\qquad\forall\beta>\lambda_1.
\eeq
Arguing again by contradiction, assume that (up to a subsequence) 
$\lim_{\e\to 0}\|u^+_\be\|_{L^2(\Omega)}=+\infty$ for some
$\b>\lambda_1$.
Then, as $\e\to 0$, from (\ref{eq}) we obtain
\beq
\io [D\bar u^+_\b D\psi -\b \bar
u^+_\b\psi]dx=0\qquad\forall\psi\in\H.
\eeq
Thus, we again have a contradiction because $D\bar u^+_\b\not\equiv 0$
on $\Omega\cap\partial (\supp \bar u_\b^+)$.

\no Therefore, we can say that for all $\b>\lambda_1$ (up to a
subsequence) $u_\be$ converges as $\e\to 0$ to a sign changing
function $u_\b$ strongly in $\H$.
Moreover, if we let $\e\to 0$ in (\ref{eq}), because of the minimality
property (\ref{min}) given by Proposition \ref{P2.3}, we infer that
$u_\b$ is a weak solution of the Dirichlet problem
\beq
\Delta u_\b-\alpha_\beta u^-_\b+\b u^+_\b=0\quad\mbox{ in
}\Omega,\qquad u_\b=0\quad\mbox{ on }\partial\Omega.
\eeq
It follows that $\io |Du^+_\b|^2dx=\b\io (u^+_\b)^2dx>0$ and 
$\io |Du^-_\b|^2dx=\alpha_\b$ $\forall\b>\lambda_1$ (notice that 
$\io |Du^-_\b|^2dx>\lambda_1$ because $u_\b^+\not\equiv 0$).
So we can say that the function $\bar
u_\b=-u_\b^-+\|u^+_\b\|_{L^2(\Omega)}^{-1}u_\b^+$ is a minimizing
function for the infimum in (\ref{ab}) and $\alpha_\beta\le\alpha$ for
every $\alpha>\lambda_1$ such that $(\alpha,\beta)\in\Sigma$, which
completes the proof. 

\qed

\begin{prop}
\label{P2.5}
For all $\b>\lambda_1$, let $\alpha_\b$ be the positive number
introduced in Theorem \ref{T2.1}.
Then $\alpha_\beta$ is continuous and strictly decreasing with respect
to $\beta$ in $]\lambda_1,+\infty[$.
Moreover, $\lim_{\b\to +\infty}\alpha_\b=\lambda_1$ and 
$\lim_{\b\to \lambda_1}\alpha_\b=+\infty$.
\end{prop}

\proof
First, let us prove that $\alpha_\b$ depends continuously on $\b$ in
$]\lambda_1,+\infty[$, namely, $\lim_{\b\to
  \bar\b}\alpha_\b=\alpha_{\bar\b}$ $\forall\bar \beta\in
]\lambda_1,+\infty[$.

\no Let us set $\bar
u_\beta=-u^-_{\b}+\|u^+_{\b}\|^{-1}_{L^2(\Omega)}u^+_{\b}$
$\forall\beta>\lambda_1$, where $u_\beta$ is the eigenfunction,
corresponding to the pair $(\alpha_\beta,\beta)$, given by Proposition
\ref{P2.4}.

\no In order to prove that
$\liminf_{\b\to\bar\b}\alpha_\b\ge\alpha_{\bar\b}$, we argue by
contradiction and assume that there exists a sequence $(\beta_n)_n$
such that $\lim_{n\to\infty}\b_n=\bar\b$ and
$\lim_{n\to\infty}\alpha_{\b_n}<\alpha_{\bar\b}$. 

\no Since $\alpha_{\b_n}=\io|D\bar u^-_{\b_n}|^2dx$ and $\io|D\bar
u^+_{\b_n}|^2dx  =\b_n$ $\forall n\in\N$, it follows that the sequence
$(\bar u_{\b_n})_n$ is bounded in $H_0^1(\Omega)$, so (up to a
subsequence) it converges to a function $\hat u$ in $L^2(\Omega)$,
weakly in $H^1_0(\Omega)$ and a.e. in $\Omega$.
Then, we have
\beq
\lim_{n\to\infty}\io |D\bar u_{\b_n}^-|^2dx\ge \io |D\hat u^-|^2dx
\quad \mbox{ and} \quad 
 \io |D\hat u^+|^2dx\le \lim_{n\to\infty}\io |D\bar
 u_{\b_n}^+|^2dx=\bar \b.
\eeq
As a consequence, there exists a function $\hat v\in H_0^1(\Omega)$
such that $0\le \hat v\le\hat u^+$ in $\Omega$ and $\io |D\hat
v|^2dx=\bar\b\io \hat v^2dx>0$.

\no Therefore, by (\ref{ab}) we have 
\beq
\alpha_{\bar\b}\le\io |D\hat u^-|^2dx\le\lim_{n\to\infty}\io |D\bar
  u_{\b_n}^-|^2dx=\lim_{n\to\infty}\alpha_{\b_n},
\eeq
which is a contradiction. 
Thus we can say that
$\liminf_{\b\to\bar\b}\alpha_\b\ge\alpha_{\bar\b}$.

\no Now, let us prove that $\limsup_{\b\to\bar\b}
\alpha_\b\le\alpha_{\bar\b}$.
By slight perturbations of the function $\bar u_{\bar\b}$, one can
construct $\forall\b>\lambda_1$ a function $\tilde u_\b\in
H^1_0(\Omega)$ such that $\|\tilde u_\b^+\|_{L^2(\Omega)}=\|\tilde
u_\b^-\|_{L^2(\Omega)}=1$, $\io |D\tilde u_\b^+|^2dx=\b$
$\forall\beta>\lambda_1$ and $\tilde u^-_\beta\to\bar u^-_{\bar\b}$ in
$H^1_0(\Omega)$ as $\b\to\bar\b$.
It follows that $\alpha_\b\le\io |D\tilde u_\b^-|^2dx$
$\forall\b>\lambda_1$ and 
\beq
\limsup_{\b\to\bar\b}\alpha_\b\leq\lim_{\b\to\bar\b}\io |D\tilde
u_\b^-|^2dx=\io |D\bar u^-_{\bar\b}|^2dx=\alpha_{\bar\b}.
\eeq
Thus, we can conclude that $\alpha_\b$ depends continuously on $\b$ in
$]\lambda_1,+\infty[$.

\no Now, we prove that $\alpha_\b$ is strictly decreasing with respect
to $\b$ in $]\lambda_1,+\infty[$.
Let us consider $\b_1$ and $\b_2$ in $]\lambda_1,+\infty[$ such that
$\b_1<\b_2$.
Since $\alpha_{\b_1}=\io |Du^-_{\b_1}|^2dx$ and $\b_1<\b_2$, there
exists $\check u_{\b_2}\in\H$ (with $\supp(\check u_{\b_2}^+)\subseteq
\supp( u_{\b_1}^+)$ and $\supp(\check u_{\b_2}^-)\supseteq \supp( u_{\b_1}^-)$)
such that $\|\check u_{\b_2}^+\|_{L^2(\Omega)}=\|\check
u_{\b_2}^-\|_{L^2(\Omega)}=1$, $\io |D \check u_{\b_2}^+|^2dx=\b_2$,
$\io|D \check u_{\b_2}^-|^2dx<\io|D  u_{\b_1}^-|^2dx=\alpha_{\b_1}$.

\no Therefore, taking into account Proposition \ref{P2.3} (see
(\ref{min})), we obtain
\begin{eqnarray}
\alpha_{\b_2}& = \min &\left\{\io |Du^-|^2dx\ :\ u\in\H,\
  \|u^+\|_{L^2(\Omega)}= \|u^-\|_{L^2(\Omega)}=1,\right.
 \nonumber \\
& & \hspace{1cm} \left. \io
  |Du^+|^2dx=\b_2\right\}< \io |Du^-_{\b_1}|^2dx
\end{eqnarray}
that is $\alpha_{\b_2}<\alpha_{\b_1}$.

\no Let us prove that $\lim_{\b\to +\infty}\alpha_\b=\lambda_1$.
Since $\alpha_\b$ is decreasing with respect to $\b$, from Proposition
\ref{P2.3} we obtain
$$
\lambda_1< \alpha_\b =  \inf \left\{\io |Du^-|^2dx\ :\ u\in\H,\
  \|u^+\|_{L^2(\Omega)}= \|u^-\|_{L^2(\Omega)}=1,\right.
$$
\beq
\hspace{2cm}\left. \io
  |Du^+|^2dx\le\b\right\}.
\eeq
As $\b\to +\infty$, it follows that
\begin{eqnarray}
\lambda_1 & \le & \lim_{\b\to +\infty}\alpha_\b =  \inf\{\alpha_\beta\ :\ \b>\lambda_1\}
\nonumber \\
& \le &
\inf \left\{\io |Du^-|^2dx\ :\ u\in\H,\
  \|u^+\|_{L^2(\Omega)}= \|u^-\|_{L^2(\Omega)}=1\right\}
\nonumber \\
& \le &
\int_{\Omega_r}|D e_1(\Omega_r)|^2dx\qquad \forall r\in\left]0,{1\over
    2}\diam(\Omega)\right[,
\end{eqnarray}
where $\Omega_r$ and $e_1(\Omega_r)$ are defined in the following way:
we fix $x_0\in\partial\Omega$, set $\Omega_r=\Omega\setminus\overline
B(x_0,r)$ and denote by $e_1(\Omega_r)$ the positive eigenfunction,
normalized in $L^2(\Omega_r)$, corresponding to the first eigenvalue of
the operator $-\Delta$ in $H^1_0(\Omega_r)$.
Then, as $r\to 0$, we obtain $\lim_{\b\to
  +\infty}\alpha_\b=\lambda_1$.
Notice that, as a consequence, we have also
$\lim_{\b\to\lambda_1}\alpha_\beta=+\infty$, because the set
$\{(\alpha_\b,\b)\in\R^2$ : $\b>\lambda_1\}$ (the first nontrivial
curve of $\Sigma$) is symmetric with respect
to the line $\{(\alpha,\beta)\in\R^2$ : $\alpha=\beta\}$ (since a pair
$(\alpha,\beta)\in\Sigma$ if and only if $(\beta,\alpha)\in\Sigma$).

\qed

{\underline {\sf Proof of Theorem \ref{T2.1}}} \hspace{2mm}
For all $\b>\lambda_1$ and $\e>0$, let $u_\be$ be a minimizing
function for the functional $f_\be$ on the set $M^\be$ (here we use
Proposition \ref{P2.2}).
From Propositions \ref{P2.3} and \ref{P2.4} we deduce that, as $\e\to
0$, $u_\be$ converges in $\H$ to a function $u_\b$ such that
$u^+_\b\not\equiv 0$ and $\|u^-_\b\|_{L^2(\Omega)}=1$, satisfying the
equation $\Delta u_\b-\alpha_\b u^-_\b+\b u^+_\b=0$ in $\Omega$, with
$\alpha_\b=\io |Du^-_\b|^2dx>\lambda_1$.
Thus $(\alpha_\b,\b)\in\Sigma$.
Moreover, from Proposition \ref{P2.3} we infer that the function $\bar
u_\b=-u^-_\b+\|u^+_\b\|^{-1}_{L^2(\Omega)}u^+_\b$ is a minimizing
function for (\ref{ab}) and  $\alpha_\b\le\alpha$
for every $\alpha>\lambda_1$ such that $(\alpha,\b)\in\Sigma$.
Notice that the eigenfunction $u_\b$ corresponding to the pair
$(\alpha_\b,\b)$ may be written as $u_\b=-\bar u_\b^-+\mu_\b\bar
u_\b^+$ with $\mu_\b=\|u^+_\b\|_{L^2(\Omega)}$.

\no By Proposition \ref{P2.5}, we know that $\alpha_\b$ is continuous
and strictly decreasing with respect  to $\b$ in $]\lambda_1,+\infty[$
and that $\lim_{\b\to+\infty}\alpha_\b=\lambda_1$,
$\lim_{\b\to\lambda_1}\alpha_\b=+\infty$.
As a consequence, we can infer that $\alpha_{\lambda_2}=\lambda_2$.
In fact, because of the minimality property (\ref{ab}), since
$(\lambda_2,\lambda_2)\in\Sigma$, we have
$\alpha_{\lambda_2}\le\lambda_2$.
Arguing by contradiction, assume that $\alpha_{\lambda_2}<\lambda_2$.
Then, the continuous curve $\{(\alpha_\b,\b)\in\R^2$ :
$\b\in]\lambda_1,+\infty[\}$ meets the line $\{(\alpha,\b)\in\R^2$ :
$\alpha=\beta\}$ in a point $(\bar\lambda,\bar\lambda)\in\R^2$ with
$\lambda_1<\bar\lambda<\lambda_2$, which is impossible because
$\bar\lambda$ must be an eigenvalue for the Laplace operator $-\Delta$
in $\H$.
Thus, we can conclude that $\alpha_{\lambda_2}=\lambda_2$.

\no Now let us prove that, as $\b\to +\infty$, $u^-_\b\to e_1$ in
$\H$.
In fact, $\|u^-_\b\|_{L^2(\Omega)}=1$ and $\alpha_\b=\io
|Du^-_\b|^2dx$ $\forall \b>\lambda_1$.
Since $\lim_{\b\to +\infty}\alpha_\b=\lambda_1$, $u^-_\b$ converges to
a function $\bar u\in \H$ in $L^2(\Omega)$, weakly in $\H$ and a.e. in
$\Omega$; moreover,
\beq
\lambda_1=\lim_{\b\to +\infty}\io |D u^-_\b|^2dx\le\io|D\bar u|^2dx
\eeq
because $\bar u\in \H$ and $\|\bar u\|_{L^2(\Omega)}=1$.
It follows that $u^-_\b\to\bar u$ strongly in $\H$ and $\bar u=e_1$.
Finally, taking into account that a function $u$ satisfies (\ref{F}) if
and only if $-u$ satisfies (\ref{F}) with $(\b,\alpha)$ in place of
$(\alpha,\b)$, we infer that $\|u^+_\b\|^{-1}_{L^2(\Omega)}u^+_\b\to
e_1$ in $\H$, as $\b\to\lambda_1$, so the proof is complete.

\qed


\sezione{Comparison with other curves of $\Sigma$ and final remarks}


Now, our aim is to use the variational characterization of the first
nontrivial curve of the Fu\v{c}\'{\i}k spectrum $\Sigma$, obtained
in Section 2, in order to prove that this curve is distinct from all
the infinitely many curves of $\Sigma$ obtained in previous papers
(see \cite{c,im} and the references therein).

\no In the next theorem we gather the main results presented in \cite{c}
and \cite{im}.

\begin{teo}
\label{T3.1}
Let $\Omega$ be a bounded connected domain of $\R^N$ with $N\ge 2$.
Then, there exists a nondecreasing sequence $(b_k)_k$ of positive
numbers, having the following properties.
For every positive integer $k$ and for all $\b>b_k$, there exist
$\alpha_{k,\b}>\lambda_1$ and $u_{k,\b}\in\H$, with
$u^+_{k,\b}\not\equiv 0$ and $u^-_{k,\b}\not\equiv 0$, such that the
equation $\Delta u_{k,\b}-\alpha_{k,\b}u^-_{k,\b}+\b u^+_{k,\b}=0$ in
$\Omega$ is satisfied for all $\b>b_k$.
Moreover, for every positive integer $k$, $\alpha_{k,\b}$ depends
continuously on $\b$ in $]b_k,+\infty[$,
$\alpha_{k,\b}<\alpha_{k+1,\b}$ $\forall\b>b_{k+1}$,
$\alpha_{k,\b}\to\lambda_1$ as $\b\to +\infty$, while $u_{k,\b}\to
-e_1$ in $\H$.
\end{teo}

\no Thus, the continuous curves $\{(\alpha_{k,\b},\b)\in\R^2$ : $\b>b_k\}$
and $\{(\b,\alpha_{k,\b})\in\R^2$ : $\b>b_k\}$ are included in the nontrivial
part of the Fu\v{c}\'{\i}k spectrum $\Sigma$ for all $k\in\N$ and, as
the first nontrivial curve of $\Sigma$, are all asymptotic to the
lines $\{\lambda_1\}\times\R$ and $\R\times\{\lambda_1\}$.

\no In addition, these curves and the corresponding eigenfunctions
have the properties described in the following proposition
(see \cite{c,im} and the references therein).

\begin{prop}
\label{P3.2}
Let $\Omega$ be a bounded connected domain of $\R^N$ with $N\ge 2$.
For every positive integer $k$, let $b_k>0$ and, for $\b>b_k$, let
$\alpha_{k,\b}>\lambda_1$ and $u_{k,\b}\in\H$ be the positive number
and the function given by Theorem \ref{T3.1}.
Then, the following properties hold.
There exist $r>0$ and, for all $k\ge 1$ and $\b>b_k$, $k$ points
$x_{1,\b},\ldots,x_{k,\b}$ in $\Omega$ such that the balls
$B\left(x_{1,\b},{r\over\sqrt\b}\right)$, \ldots, 
$B\left(x_{k,\b},{r\over\sqrt\b}\right)$ are pairwise disjoint and
all included in $\Omega$, $u_{k,\b}(x)\le 0$ $\forall
x\in\Omega\setminus\cup_{i=1}^k
B\left(x_{i,\b},{r\over\sqrt\b}\right)$ and $u^+_{k,\b}\not\equiv 0$
in $B\left(x_{i,\b},{r\over\sqrt\b}\right)$ $\forall
i\in\{1,\ldots,k\}$.
As $\b\to+\infty$, we have
\beq
\lim_{\b\to +\infty}e_1(x_{i,\b})=\max_\Omega e_1\qquad\forall
i\in\{1,\ldots,k\}
\eeq
and
\beq
\lim_{\b\to +\infty}\sqrt\b\, |x_{i,\b}-x_{j,\b}|=\infty\quad\mbox{ for
}i\neq j.
\eeq
If, $\forall k\in\N$, $\forall i\in\{1,\ldots,k\}$, $\forall\b>b_k$,
we set $s_{i,k,\b}=\sup\left\{u_{k,\b}(x)\ :\ x\in
B\left(x_{i,\b},{r\over\sqrt\b}\right)\right\}$ and define
$U_{i,k,\b}(x)=s^{-1}_{i,k,\b}u_{k,\b}\left({x\over\sqrt\b}+x_{i,\b}\right)$
$\forall x\in \sqrt\b(\Omega-x_{i,\b})$, then the rescaled function
$U_{i,k,\b}$ converges as $\b\to +\infty$ to the radial solution $U$
of the problem 
\beq
\Delta U+U^+=0\quad\mbox{ in }\R^N,\qquad U(0)=1
\eeq
and the convergence is uniform on the compact subsets of $\R^N$.

\no If $N\ge 3$, we have 
\beq
\label{kb3}
\lim_{\b\to+\infty}\b^{N-2\over 2}(\alpha_{k,\b}-\lambda_1)=\cp(\bar
r_1)\left(\max_\Omega e_1\right)^2k
\eeq
where $\bar r_1$ is the radius of the balls in $\R^N$ for which the
first eigenvalue of $-\Delta$ in $H^1_0$ is equal to 1 and $\cp(\bar
r_1)$ denotes the capacity of these balls.

\no Finally, in the case $N=2$ we have
\beq
\label{kb2}
\lim_{\b\to +\infty} \lg\b(\alpha_{k,\b}-\lambda_1)=4\pi
\left(\max_\Omega e_1\right)^2k.
\eeq
\end{prop}

\no In Proposition \ref{P2.5} we proved that $\alpha_\b$ converges to
$\lambda_1$ as $\b\to +\infty$; now, we need to estimate the rate of
convergence. 

\begin{prop}
\label{P3.3}
For all $\b>\lambda_1$, let $\alpha_\b$ be the positive number
introduced in Theorem \ref{T2.1}.
Then, for $N\ge 3$ we have
\beq
\label{b3}
\lim_{\b\to +\infty}\b^{N-2\over 2}(\alpha_\b-\lambda_1)=0
\eeq
while, for $N=2$,
\beq
\label{b2}
\lim_{\b\to +\infty}\lg\b (\alpha_\b-\lambda_1)=0.
\eeq
\end{prop}

\proof
Let $\bar r_1$ be the positive number introduced in Proposition
\ref{P3.2}. 
For all $y\in\Omega$, let us consider the function $\bar
u_{\b,y}\in\H$, defined as follows.
First notice that $B\left(y,{\bar
    r_1\over\sqrt\b}\right)\subseteq\Omega$ for $\b$ large enough
and the first eigenvalue
of $-\Delta$ in $H^1_0\left(B\left(y,{\bar
      r_1\over\sqrt\b}\right)\right)$ is equal to $\b$ (because of the
choice of $\bar r_1$).
Then, for $\b>0$ large enough, in the ball $B\left(y,{\bar
    r_1\over\sqrt\b}\right)$ we define $\bar u_{\b,y}$ to be the
positive eigenfunction corresponding to the first eigenvalue of
$-\Delta$ in 
$H^1_0\left(B\left(y,{\bar r_1\over\sqrt\b}\right)\right)$, 
normalized in $L^2\left(B\left(y,{\bar r_1\over\sqrt\b}\right)\right)$.
Now, in order to define $\bar u_{\b,y}$ in $\Omega \setminus
B\left(y,{\bar r_1\over\sqrt\b}\right)$, set $\e_\b=\b^{-q}$ with
$q\in\left]{1\over 2}-{1\over N},{1\over 2}\right[$ if $N\ge 3$ and
$\e_\b={1\over \lg\b}$ if $N=2$.
Then, for $\b>0$ large enough, consider the function $\tilde u_{\b,y}$
in $H^1_0\left(\Omega \setminus
\overline B\left(y,{\bar r_1\over\sqrt\b}\right)\right)$, such that
$\Delta \tilde 
u_{\b,y}=0$ in the annulus $A\left(y,{\bar r_1\over \sqrt\b},{\bar
    r_1\over \sqrt\b}+\e_\b\right)= B\left(y,{\bar r_1\over
    \sqrt\b}+\e_\b\right)\setminus\overline B \left(y,{\bar r_1\over
    \sqrt\b}\right)$ and $\tilde u_{\b,y}(x)=-e_1(x)$ $\forall x\in
\Omega\setminus B\left(y,{\bar r_1\over
    \sqrt\b}+\e_\b\right)$.
Finally, we complete the definition of $\bar u_{\b,y}$ by setting
$\bar u_{\b,y}(x)=\|\tilde
u_{\b,y}\|^{-1}_{L^2(\Omega\setminus\overline B(y,\bar
  r_1/\sqrt\b))}\tilde u_{\b,y}(x)$ $\forall x\in
\Omega\setminus\overline B\left(y,{\bar r_1\over\sqrt\b}\right)$.

\no Since  $\|\bar u^+_{\b,y}\|_{L^2(\Omega)}=\|\bar
u^-_{\b,y}\|_{L^2(\Omega)}=1$ and $\io |D\bar u^+_{\b,y}|^2dx=\b$,
taking into account (\ref{ab}) we infer
that $\alpha_\beta\le\io |D\bar u^-_{\b,y}|^2dx$ for $\b>0$ large enough
so that $B\left(y,{\bar r_1\over
    \sqrt\b}+\e_\b\right)\subseteq\Omega$.

\no Let us estimate the integral $\io |D\bar u^-_{\b,y}|^2dx$.
We have
\beq
 \io |D\bar u^-_{\b,y}|^2dx=\|\tilde
 u_{\b,y}\|^{-2}_{L^2(\Omega\setminus\overline B(y,\bar r_1/\sqrt\b))}
\int_{\Omega\setminus\overline B\left(y,{\bar r_1\over \sqrt\b}\right)}
|D \tilde u_{\b,y}|^2dx
\eeq
where
\beq
\label{L2}
\|\tilde u_{\b,y}\|^{2}_{L^2(\Omega\setminus\overline B(y,\bar r_1/\sqrt\b))}
=1-\int_{B\left(y,{\bar r_1\over\sqrt\b}+\e_\b\right)}e_1^2dx+
\int_{A\left(y,{\bar r_1\over \sqrt\b},{\bar r_1\over \sqrt\b}+\e_\b\right)}
\tilde u^2_{\b,y}dx
\eeq
and
\beq
\label{H1}
\int_{\Omega\setminus\overline B\left(y,{\bar r_1\over\sqrt\b}\right)} |D \tilde
u_{\b,y}|^2dx=
\lambda_1-
\int_{B\left(y,{\bar r_1\over\sqrt\b}+\e_\b\right)}|D e_1|^2dx+
\int_{A\left(y,{\bar r_1\over \sqrt\b},{\bar r_1\over \sqrt\b}+\e_\b\right)}
|D\tilde u_{\b,y}|^2dx.
\eeq
Since $\lim_{\b\to +\infty} \e_\b\sqrt\b=+\infty$, one can easily
verify that there exists a positive number $c(y)$ (depending only on
$y$) such that
\beq
\label{cy1}
\lim_{\b\to+\infty}\e^{-N}_\b\int_{B\left(y,{\bar r_1\over\sqrt\b}+\e_\b\right)}
(e_1^2+|D e_1|^2)dx\le c(y)
\eeq
and
\beq
\label{cy2}
\limsup_{\b\to +\infty}\e_\b^{-N}
\int_{A\left(y,{\bar r_1\over \sqrt\b},{\bar r_1\over
      \sqrt\b}+\e_\b\right)}
\tilde u_{\b,y}^2dx\le c(y).
\eeq
Now, let us estimate the integral $\int_{A\left(y,{\bar r_1\over
      \sqrt\b},{\bar r_1\over \sqrt\b}+\e_\b\right)}|D\tilde
u_{\b,y}|^2dx$.
Let us write $\tilde u_{\b,y}$ as $\tilde u_{\b,y}=\tilde
v_{\b,y}+\tilde w_{\b,y}$ where $\tilde v_{\b,y}$ and $\tilde
w_{\b,y}$ satisfy $\Delta \tilde v_{\b,y}=0$,  $\Delta \tilde w_{\b,y}=0$
in $A\left(y,{\bar r_1\over\sqrt\b},{\bar r_1\over \sqrt\b}+\e_\b\right)$
with boundary condition $\tilde v_{\b,y}=e_1$,  $\tilde w_{\b,y}=-e_1$
on $\partial B\left(y,{\bar r_1\over\sqrt\b}\right)$ and  $\tilde
v_{\b,y}=0$,  $\tilde w_{\b,y}=-e_1$ on $\partial B\left(y,{\bar
    r_1\over\sqrt\b}+\e_\b\right)$.

\no If $N\ge 3$, one can verify by standard arguments that
\beq
\lim_{\b\to +\infty}\b^{N-2\over 2}
\int_{A\left(y,{\bar r_1\over \sqrt\b},{\bar r_1\over \sqrt\b}+\e_\b\right)}
|D\tilde v_{\b,y}|^2dx= e_1^2(y)\int_{\R^N\setminus\overline B(0,\bar
  r_1)}|D\tilde V|^2dx=e_1^2(y)\cp(\bar r_1),
\eeq
where $\tilde V$ satisfies $\Delta \tilde V=0$ in
$\R^N\setminus\overline B(0,\bar r_1)$, $\tilde V=1$ on $\partial
B(0,\bar r_1)$, $\lim_{|x|\to\infty}\tilde V(x)=0$.

\no Moreover,
\beq
\limsup_{\b\to +\infty} \e^{-N}_\b
\int_{A\left(y,{\bar r_1\over \sqrt\b},{\bar r_1\over
      \sqrt\b}+\e_\b\right)}
|D\tilde w_{\b,y}|^2 dx\le c(y),
\eeq
where $c(y)$ is the positive number appearing in (\ref{cy1}) and
(\ref{cy2}).
Taking into account the choice of $\e_\b$, we obtain also
\beq
\lim_{\b\to +\infty}\b^{N-2\over 2}
\int_{A\left(y,{\bar r_1\over \sqrt\b},{\bar r_1\over
      \sqrt\b}+\e_\b\right)}
|D\tilde v_{\b,y}|\,|D\tilde w_{\b,y}|\,dx=0.
\eeq
Thus, since $\alpha_\b\le\io |D\bar u_{\b,y}^-|^2dx$ for $\b>0$ large
enough so that $B\left(y,{\bar r_1\over \sqrt\b}+\e_\b\right)\subset\Omega$, it follows that
\beq
\limsup_{\b\to +\infty}\b^{N-2\over 2}(\alpha_\b-\lambda_1)\le
e_1^2(y)\cp(\bar r_1)\qquad\forall y\in\Omega,
\eeq
which (as $y$ tends to the boundary of $\Omega$) implies (\ref{b3}).

\no If $N=2$, we argue in analogous way (but with $\lg \b$ in place
of $\b^{N-2\over 2}$).
Since $\e_\b={1\over \lg \b}$ for $N=2$, one can verify by direct
computation that
\beq
\lim_{\b\to +\infty}\lg\b
\int_{A\left(y,{\bar r_1\over \sqrt\b},{\bar r_1\over
      \sqrt\b}+\e_\b\right)}
|D\tilde v_{\b,y}|^2dx=4\pi e_1^2(y).
\eeq
Moreover,
\beq
\limsup_{\b\to +\infty}\lg^2\b
\int_{A\left(y,{\bar r_1\over \sqrt\b},{\bar r_1\over
      \sqrt\b}+\e_\b\right)}
|D\tilde w_{\b,y}|^2dx\le c(y)
\eeq
(where $c(y)$ is the same as in (\ref{cy1}) and (\ref{cy2})) and
\beq
\lim_{\b\to +\infty}\lg\b
\int_{A\left(y,{\bar r_1\over \sqrt\b},{\bar r_1\over
      \sqrt\b}+\e_\b\right)}
|D\tilde v_{\b,y}|\,|D\tilde w_{\b,y}|\,dx=0.
\eeq
Therefore, it follows that
\beq
\limsup_{\b\to +\infty}\lg \b(\alpha_\b-\lambda_1)\le
4\pi e_1^2(y)\qquad\forall y\in\Omega,
\eeq
which clearly implies (\ref{b2}).

\qed

As a consequence of Propositions \ref{P3.2} and \ref{P3.3} we can
state the following corollary.

\begin{cor}
\label{C3.4}
For every positive integer $k$, there exists $\tilde b_k>0$ such that
$\alpha_\b<\alpha_{k,\b}$ $\forall \b>\tilde b_k$.
\end{cor}

\no The proof follows directly by comparing formulas (\ref{kb3}) and
(\ref{kb2}) with (\ref{b3}) and (\ref{b2}).

\begin{rem}
\label{R3.5}
{\em
The proof of Proposition \ref{P3.3} suggests that, for $\b>0$ large
enough, the support of $u^+_\b$ is localized near the boundary of
$\Omega$. 
Indeed, a more careful analysis of the asymptotic behaviour of $u_\b$
as $\b\to +\infty$ (arguing as in \cite{im}) shows that for all $\b>0$
there exists $y_\b\in\Omega$ such that (up to a subsequence) $y_\b$
converges, as $\b\to +\infty$, to a point $\bar y\in\partial\Omega$
and the rescaled function $\left(\sup_\Omega
  u_\b\right)^{-1}u_\b\left(y_\b+{x\over \sqrt\b}\right)$ converges
to a function $\bar U$ such that 
\beq 
\Delta\bar U+\bar U^+=0, \qquad \bar
U^+\not\equiv 0,\  \bar U^-\not\equiv 0 \quad\mbox{in } H,\qquad \bar U=0
\quad\mbox{on }\partial H
\eeq
with $H=\{x\in\R^N$ : $(x\cdot\bar\nu)<0\}$, where $\bar\nu$ denotes
the outward normal to $\partial\Omega$ in $\bar y$.

\no Moreover (arguing as in \cite{im}) for $\b$ large enough one can
construct multibumps eigenfunctions for the Fu\v{c}\'{\i}k spectrum,
having an arbitrarily large number of bumps localized near prescribed
connected components of $\partial\Omega$.

\no In fact, if $\Omega$ is a smooth bounded connected domain of
$\R^N$ with $N\ge 2$, arguing as in \cite{im}, $\forall k\in\N$ one
can construct, for $\b>0$ large enough, a $k$-bumps eigenfunction
$\bar u_{k,\b}$, corresponding to a pair
$(\bar\alpha_{k,\b},\b)\in\Sigma$, with $k$ bumps localized near
$\partial \Omega$ and having asymptotic profile described by the
functions $\bar U$.
Here $\bar\alpha_{k,\b}$ depends
continuously on $\b$ and, for $\b>0$ large enough, we have
$\bar\alpha_{k,\b}<\bar\alpha_{k+1,\b}$ $\forall k\in\N$ and
$\bar\alpha_{k_1,\b}<\alpha_{k_2,\b}$ $\forall k_1,k_2\in\N$.
Thus, if $N\ge 2$, we obtain a new class of infinitely many curves in
$\Sigma$, asymptotic to the lines $\{\lambda_1\}\times\R$ and
$\R\times\{\lambda_1\}$, while in the case $N=1$ there exist only two
curves having this property.

}\end{rem}

\begin{rem}
{\em
 Notice that the difference between the cases $N=1$ and $N\ge 2$ is
even more evident if we replace the Dirichlet boundary condition by
the Neumann condition ${\partial u\over\partial\nu}=0$ on
$\partial\Omega$.
In fact, if we denote by
$\tilde\lambda_1<\tilde\lambda_2\le\tilde\lambda_3\le\ldots$ and by
$\widetilde\Sigma$, respectively, the eigenvalues of the Laplace
operator $-\Delta$ and the Fu\v{c}\'{\i}k spectrum  with Neumann
boundary conditions, we have $\tilde\lambda_1=0$ and, if $N=1$, no
curve of $\widetilde\Sigma$ is asymptotic to the lines
$\{0\}\times\R$ or $\R\times\{0\}$.
Indeed, a direct computation shows that, for $N=1$, every nontrivial pair
$(\alpha,\beta)$ of $\widetilde\Sigma$ satisfies
$\alpha>{1\over 4}\tilde\lambda_2$ and $\b>{1\over 4}\tilde\lambda_2$,
with $\tilde\lambda_2>0$.
On the contrary, in the case $N\ge 2$ there exist infinitely many
curves contained in $\widetilde\Sigma$  and asymptotic to the lines 
$\{0\}\times\R$ or $\R\times\{0\}$.
The corresponding eigenfunctions have an arbitrarily large number of
bumps localized in the interior of $\Omega$ or near prescribed
connected components of $\partial\Omega$.
Both, interior and boundary bumps, present the same asymptotic
profile, described by the function $U$ introduced in Proposition
\ref{P3.2} (while, in case of Dirichlet boundary conditions, the
asymptotic profile is described by the function $U$ for the interior
bumps and by the function $\bar U$ for the boundary bumps).

\no Notice that, as pointed out in \cite{dFG,CdFG,ACG}, in case of Neumann
boundary conditions there exists a strict connection between the
nonexistence of curves in  $\widetilde\Sigma$, asymptotic to the
lines $\{0\}\times\R$ and $\R\times\{0\}$, and the fact that the
antimaximum principle (see \cite{CP}) holds uniformly (in a suitable
sense).
}
\end{rem}

\begin{rem}{\em
\label{R3.6}
For the sake of simplicity, in this paper we have considered only the case of
the Laplace operator, but the results we have presented may be easily
extended to cover the case of more general boundary conditions and
elliptic operators in divergence form.
Moreover, the variational method we have used in this paper may be easily
adapted to deal also with quasilinear operators as the $p$-laplacian.
Thus, we can obtain also for the $p$-laplacian a variational
characterization of the first nontrivial curve of the Fu\v{c}\'{\i}k
spectrum, similar to Theorem \ref{T2.1} (but the asymptotic behaviour
depends on $p$, on the spatial dimension $N$ and on the boundary
conditions we consider).
}\end{rem}

{\small {\bf Acknowledgement}. The authors have been supported by the ``Gruppo
Nazionale per l'Analisi Matematica, la Probabilit\`a e le loro
Applicazioni (GNAMPA)'' of the {\em Istituto Nazionale di Alta Matematica
(INdAM)}.}



\end{document}